\documentclass[12pt,a4paper]{amsart}
\usepackage{aliascnt}
\usepackage{hyperref}
\usepackage{amsfonts}
\usepackage{mathrsfs}
\usepackage{amsmath}
\usepackage{amsmath,graphics}%diagrams,
\usepackage{amssymb}
\usepackage{indentfirst,latexsym,bm,amsmath,pstricks,amssymb,amsthm,graphicx,fancyhdr,float,color}
\usepackage[all]{xy}
\usepackage{amsthm}
\usepackage{amscd}
\usepackage[all]{hypcap}

\newtheorem{theorem}{Theorem}[section]

\newaliascnt{lemma}{theorem}
\newtheorem{lemma}[lemma]{Lemma}
\aliascntresetthe{lemma}

\newaliascnt{conjecture}{theorem}

\aliascntresetthe{conjecture}

\newaliascnt{proposition}{theorem}
\newtheorem{proposition}[proposition]{Proposition}
\aliascntresetthe{proposition}

\newaliascnt{corollary}{theorem}
\newtheorem{corollary}[corollary]{Corollary}
\aliascntresetthe{corollary}

\newaliascnt{problem}{theorem}

\aliascntresetthe{problem}

\newaliascnt{claim}{theorem}

\aliascntresetthe{claim}

\theoremstyle{definition}

\newaliascnt{definition}{theorem}
\newtheorem{definition}[definition]{Definition}
\aliascntresetthe{definition}

\newaliascnt{example}{theorem}
\newtheorem{example}[example]{Example}
\aliascntresetthe{example}

\theoremstyle{remark}

\newaliascnt{remark}{theorem}
\newtheorem{remark}[remark]{Remark}
\aliascntresetthe{remark}

\newaliascnt{remarks}{theorem}

\aliascntresetthe{remarks}

\numberwithin{equation}{section}

\numberwithin{figure}{section}

\def\ol{\overline}
\def\ra{\rightarrow}
\def\lra{\longrightarrow}

\def\({$($}
\def\){$)$}

\def\rank{\text{{\rm rank\,}}}

\newcommand{\Ecal}{\mathcal{E}}
\newcommand{\Fcal}{\mathcal{F}}
\newcommand{\Mcal}{\mathcal{M}}

\newcommand{\Acal}{\mathcal{A}}

\newcommand{\Hcal}{\mathcal{H}}

\newcommand{\Lcal}{\mathcal{L}}
\newcommand{\Jac}{\mathrm{Jac}}
\newcommand{\Tcal}{\mathcal{T}}

\newcommand{\Scal}{\mathcal{S}}

\newcommand{\Wcal}{\mathcal{W}}

\newcommand{\Cbb}{\mathbb{C}}
\newcommand{\Gbb}{\mathbb{G}}

\newcommand{\Pbb}{\mathbb{P}}
\newcommand{\Qbb}{\mathbb{Q}}
\newcommand{\Rbb}{\mathbb{R}}
\newcommand{\Sbb}{\mathbb{S}}

\newcommand{\Zbb}{\mathbb{Z}}

\newcommand{\Gbf}{\mathbf{G}}
\newcommand{\Hbf}{\mathbf{H}}

\newcommand{\mrm}{\mathrm{m}}

\newcommand{\Mbar}{{\ol M}}

\newcommand{\Res}{\mathrm{Res}}
\newcommand{\bsh}{\backslash}
\newcommand{\ad}{\mathrm{ad}}
\newcommand{\der}{\mathrm{der}}
\newcommand{\isom}{\simeq}
\newcommand{\mono}{\hookrightarrow}

\newcommand{\GSp}{\mathbf{GSp}}
\newcommand{\Sp}{\mathbf{Sp}}

\newcommand{\GL}{\mathbf{GL}}

\newcommand{\SU}{\mathbf{SU}}
\newcommand{\SO}{\mathbf{SO}}

\newcommand{\tr}{\mathrm{tr}}

\newcommand{\Std}{\mathrm{Std}}

\newcommand{\Acalbar}{\overline{\Acal}}
\newcommand{\Vcal}{\mathcal{V}}

\newcommand{\Tan}{\mathrm{Tan}}
\newcommand{\maxx}{\mathrm{max}}
\newcommand{\Grbf}{\mathbf{Gr}}

\newcommand{\vbar}{\bar{v}}
\newcommand{\zbar}{\bar{z}}

\newcommand{\Mcalbar}{\overline{\Mcal}}

\newcommand{\Spin}{\mathbf{Spin}}

\newcommand{\even}{\mathrm{even}}
\newcommand{\odd}{\mathrm{odd}}

\newcommand{\Wbar}{\overline{W}}
\newcommand{\ibf}{\mathbf{i}}

\newcommand{\Pcal}{\mathcal{P}}

\title{On Higgs bundles over Shimura varieties of ball quotient type}

\author{Ke Chen}

\address{Department of Mathematics, Nanjing University, Hankou Road 22, Nanjing 210093, P. R. of China}

\email{kechen@ustc.edu.cn}

\author{Xin Lu}
\address{Institut f\"ur Mathematik, Universit\"at Mainz, Mainz, 55099, Germany}
\email{x.lu@uni-mainz.de}
\thanks{This work is supported by SFB/Transregio 45 Periods, Moduli Spaces and Arithmetic of Algebraic Varieties of the DFG (Deutsche Forschungsgemeinschaft), and by
 National Key Basic Research Program of China (Grant No. 2013CB834202) and National Natural Science Foundation of China (Grant No. 11231003, No. 11301495).}

\author{Sheng-Li Tan}
\address{Department of Mathematics, East China Normal University, Dongchuan Road 500, Shanghai 200241, P. R. of China}
\email{sltan@math.ecnu.edu.cn}

\author{Kang Zuo}
\address{Institut f\"ur Mathematik, Universit\"at Mainz, Mainz, 55099, Germany}
\email{zuok@uni-mainz.de}

\begin{document}

\dedicatory{Dedicated to Prof. Ngaiming Mok on the occasion of his sixtieth birthday}
\maketitle

\begin{abstract}
We prove the generic exclusion of certain Shimura varieties of unitary and orthogonal types from the Torelli locus.
The proof relies on a slope inequality on surface fibration due to G. Xiao,
and the main result implies that certain Shimura varieties only meet the Torelli locus in dimension zero.	
%
%We study slopes of Higgs bundles associated to symplectic representations on Shimura varieties of $\SU(n,1)$-type, and refine a previous result on the generic exclusion of Shimura subvarieties of $\SU(n,1)$-type from the Torelli locus. An application to the generic exclusion of Shimura subvarieties of orthogonal type is also discussed.
\end{abstract}

\tableofcontents

\section{Introduction}

The Coleman-Oort conjecture, cf. \cite{moonen oort survey}, predicts that when the genus $g$ is large enough, the Torelli locus $\Tcal_g$ inside the Siegel modular variety $\Acal_g$ does not contain generically any Shimura subvariety of strictly positive dimension; in other words, the open Torelli locus $\Tcal_g^\circ$ does not contain any Zariski open subvariety of an arbitrary Shimura subvariety of dimension $>0$ in $\Acal_g$. Using the Andr\'e-Oort conjecture, an unconditional proof of which has recently been given in \cite{tsimerman andre oort}, this amounts to the finiteness of CM points in $\Tcal_g^\circ$ for $g$ sufficiently large, which is also the original formulation of Coleman.

In \cite{hain locally symmetric} Hain has established the conjecture for a large class of Shimura subvarieties in $\Acal_g$ that do not contain locally symmetric divisors. In particular, it holds for Shimura subvarieties uniformized by Hermitian symmetric domains of rank at least 2. His proof makes use of rigidity property of mapping class groups, namely for an arithmetic group $\Gamma$ coming from a simple $\Qbb$-group of $\Rbb$-rank at least 2, any homomorphism from $\Gamma$ to the mapping class group $\Gamma_{g,r}^n$ is of finite image. 

% page 204 in Mok book
Note that the phenomenon of rigidity in \cite{hain locally symmetric} is also related to the metric rigidity property studied by  Ngaiming Mok in \cite{mok rigidity book}:

\begin{example}
The idea can be illustrated through the case of Hilbert modular varieties treated in \cite{de jong zhang hilbert}: assume that $g\geq 5$ and let $M\subset\Acal_g$ be the Hilbert modular subvariety  parametrizing abelian varieties of dimension $g$ with real multiplication by $O_F$, where $F$ is a totally real field of degree $g$ with an order $O_F$ in $F$. Consider an extremal situation where $M$ is not only contained in $\Tcal_g^\circ$, but actually lifts to $i:M\mono\Mcal_g$, giving rise to a surjective pull-back $i^*\Omega_{\Mcal_g}^1\ra\Omega_M^1$. Note that $\Omega^1_M$ is merely semi-positive, dual to the proper semi-negativity of $\Tan_M$ studied in \cite{mok rigidity book}, and it admits non-trivial quotients which are NOT big, which is absurd because $\Omega_{\Mcal_g}^1$, and thus $i^*\Omega^1_{\Mcal_g}$ as well, is already big.\end{example}

On the other hand, for Shimura varieties uniformized by simple Hermtian symmetric domains of rank 1, such as those associated to $\SU(n,1)$ and $\Spin(N,2)$, Hodge-theoretic techniques provide complements when the rigidity property fails. In our previous work \cite{chen lu zuo compositio}, we have proved the Coleman-Oort conjecture for a class of Shimura varieties of these types. The starting point is a numerical property of semi-stable fibration of surfaces over curves, which is translated, via the Simpson correspondence, into constraints on the fundamental group representations for the Shimura varieties of interest, and one concludes by Satake's classification of rational symplectic representations. More precisely,  the following inequality due to G. Xiao will play a crucial role in our study:

\begin{theorem}[Xiao's inequality]\label{xiao's inequality} Let $f:S\ra B$ be a non-isotrovial fibration of a smooth projective algebraic surface $S$ over a smooth projective algebraic curve $B$. Assume that $f$ is generically smooth and its fibers are semi-stable curves of genus $g\geq 2$. Then holds the inequality $$12\deg f_*\omega_{S/B}\geq (2g-2+\rank A_\maxx)\mu_\maxx$$ where $\mu_\maxx$ is the maximum of slopes of vector subbundles of $\omega_{S/B}$, and $A_\maxx\subset f_*\omega_{S/B}$ is  the maximal subbundle of $f_*\omega_{S/B}$ of slope $\mu_\maxx$.

\end{theorem}

\begin{remark}
	(1) The inequality above is slightly different from the original version of Xiao \cite{xiao slope}, where he obtains $12\deg f_*\omega_{S/B}\geq (2g-2)\mu_\maxx$ and thus $$\frac{\rank F^{1,0}}{g}\leq \frac{5}{6}+\frac{1}{6g}.$$ The proof of the version in \autoref{xiao's inequality} is exactly the same as in \cite{xiao slope}, which is not reproduced in this paper.
	
	(2) In \cite{chen lu zuo compositio} we have made use of the following refined form
	$$12\deg f_*\omega_{S/B}\geq\big(4g-4-2\rank F^{1,0}\big)\mu_\maxx, \mathrm{\quad and \quad } \frac{\rank F^{1,0}}{g}\leq\frac{4}{5}+\frac{2}{5g},$$ where a decomposition $f_*\omega_{S/B}=A^{1,0}\oplus F^{1,0}$ is assumed, with $A^{1,0}$ semi-stable and ample while $F^{1,0}$ is the maximal flat part. Such a decomposition holds when $B$ is constructed out of a Shimura curve contained generically in $\Tcal_g$ due to properties of Higgs bundles on Shimura curves: only one single non-zero slope appears in the Harder-Narasimhan filtration of $f_*\omega_{S/B}$. The Shimura varieties studied in \cite{chen lu zuo compositio} do contain Shimura curves, while in the more general case of the present paper we can only resort to the inequality as in \autoref{xiao's inequality} because several different slopes might be involved.
\end{remark}

%Write $f_*\omega_{S/B}=A^{1,0}\oplus F^{1,0}$ for the decomposition of the Hodge bundle $f_*\omega_{S/B}$ into the maximal flat subbundle $F^{1,0}$ and an ample subbundle $A^{1,0}$. Then $$\mu_\maxx\geq \mu(A^{1,0})=\frac{\deg f_*\omega_{S/B}}{g-\rank F^{1,0}} \textrm{\ and\ }     \frac{\rank{F^{1,0}}}{g}\leq \frac{5}{6}+\frac{1}{6g}.$$ If moreover $A^{1,0}$ is semi-stable, then hold further $$12\deg f_*\omega_{S/B}\geq (4g-4-2\rank F^{1,0})\mu_\maxx \textrm{\ and }\frac{\rank{F^{1,0}}}{g}\leq \frac{4}{5}+\frac{2}{5g}.$$

%The main results of \cite{chen lu zuo compositio} have only made use of the rank of the flat part in the Hodge bundles for certain Shimura varieties containing Shimura curves. In the preset paper we refine the estimation using slopes of the ample part, and establish criteria generically excluding Shimura varieties of $\SU(n,1)$-type, which do not necessarily contain Shimura curves.

\begin{theorem}[Shimura varieties of $\SU(n,1)$-type]\label{main theorem shimura varieties of su(n,1)-type} Let $M\subset\Acal_V$ be a Shimura subvariety of $\SU(n,1)$-type defined over a totally real field $F$ of degree $d$, such that the corresponding symplectic representation is primary of type $\Lambda_m$ for some integer $m\in[1,n]$ of multiplicity $r$. Then $M$ is NOT contained generically in $\Tcal_g$ as long as the following inequality holds:$$\frac{n+m-1}{n}\left(\frac{n+1}{m}\cdot d+\frac{1}{2}-\frac{1}{r\binom{n}{m-1}}\right)> 12.$$ %$$\frac{n+m-1}{n}(\frac{n+1}{m}(d+1)-\frac{2}{r\binom{n}{m-1}})>12.$$

\end{theorem}

\begin{remark}
The theorem above actually implies that the Shimura variety $M$ under consideration   only meets the open Torelli locus $\Tcal_g^\circ$ at finitely many points: otherwise $M$ contains a curve $C$ which is generically contained in $\Tcal_g$, and the arguments in the proof (cf. Section 3) produces an inequality contradicting Xiao's estimation. Our naive bound on the flat part also treats the similar phenomenon for more general Shimura varieties, cf. \autoref{proposition naive bound on the flat part} and \autoref{corollary naive bound}.

\end{remark}

The theorem above is applied to the generic exclusion from the Torelli locus of some Shimura varieties of orthogonal type containing Shimura varieties of unitary type. In fact, if $h:W\times W\ra\Cbb$ is an Hermitian form over $\Cbb$ of signature $(n,1)$, then its real part is a quadratic form of signature $(2n,2)$, and one obtains a natural equivariant embedding of the corresponding Hermitian symmetric domains. Adding suitable  arithmetic constraints we expand this example into the following theorem:
\begin{theorem}[Shimura varieties of orthogonal types]\label{main theorem shimura varieties of orthogonal types} Let $M'\subset\Acal_V$ be a Shimura subvariety of $\Spin(N,2)$-type given by some Shimura datum $(\Gbf',X';X'^+)$, namely associated to some quadratic space $(W,q)$ over a totally real field $F$ of degree $d$, of signature \begin{itemize}
\item $(N,2)$ along one fixed real embedding $\sigma=\sigma_1:F\mono\Rbb$;

\item definite along the other embeddings $\sigma_2,\cdots,\sigma_d:F\mono\Rbb$.
\end{itemize}Assume that $M'$ contains a Shimura subvariety of $\SU(n,1)$-type associated to some Hermitian space $(H,h)$ over some CM quadratic extension $E$ of $F$, such that the signature of $h$ is \begin{itemize}
\item $(n,1)$ along $\sigma$;

\item definite along the other embeddings $\sigma_2,\cdots,\sigma_d$;
\end{itemize} which fits into an orthogonal direct sum decomposition $W=U\oplus\Res_{E/F}H$ with $U$ some $F$-subspace of signature $(N-2n,0)$ along $\sigma$. If $N>2n$ and the inclusion $M'\mono\Acal_V$ is defined by a symplectic representation of primary type, then $M'$ is NOT contained generically in $\Tcal_g$ as long as $d> 3+\frac{1}{m\cdot 2^{\lfloor(N+1)/2\rfloor}}-\frac{1}{2^{n+2}}$, where $m$ is the multiplicity of the spinor representation in $W\otimes_{F,\sigma}\Rbb$ for the group $\Spin(N,2)$, the only non-compact factor in $\Gbf'^\der(\Rbb)$. %, with $m$ the multiplicity of the spinor representation in $W\otimes_{F,\sigma}\Rbb$.

In particular the inequality holds whenever $d\geq 4$; if $N> 2n+4$, then it holds for $d\geq 3$.

\end{theorem} The theorem deals with the case $N>2n$, while the case $N=2n$ is treated in detail taking care of the parity of $n$ and the primary type of symplectic representations involved, cf. \autoref{proposition restriction on the real part}.

Similar to \cite{chen lu zuo compositio}, the criteria obtained involve certain representation-theoretic parameters describing the symplectic representations defining Shimura subvarieties, and they are NOT pure bounds on the genus $g$.

The paper is organized as follows. In Section 2 we collect preliminaries on Shimura subvarieties in $\Acal_g$ and prove a naive bound on the compact factors in the Mumford-Tate groups for Shimura subvarieties contained generically in the Torelli locus. Section 3 computes the slopes of certain Higgs bundles on Shimura varieties of $\SU(n,1)$-type and proves \autoref{main theorem shimura varieties of su(n,1)-type}.  Finally in Section 4 we apply the results in Section 3 to a class of Shimura varieties of orthogonal types. 

\section{A naive bound for flat Higgs subbundles}

In this paper Shimura varieties and Shimura subvarieties in $\Acal_g$ are connected algebraic varieties over $\Cbb$, following the definitions given in \cite{chen lu zuo compositio,lu zuo}:

\begin{definition}[Shimura varieties]\label{definition shimura varieties}

A (connected) Shimura datum is of the form $(\Gbf,X;X^+)$ consisting of \begin{itemize}
\item $(\Gbf,X)$ a (pure) Shimura datum in the sense of \cite{deligne pspm};

\item $X^+$ is a connected component of $X$.
\end{itemize}

Note that $X$ is a homogeneous space under $\Gbf(\Rbb)$ of homomorphisms $\Sbb\ra\Gbf_\Rbb$ subject to certain Hodge-theoretic constraints, with $\Sbb$ the Deligne torus $\Res_{\Cbb/\Rbb}\Gbb_\mrm$; $X^+$ is homogeneous under $\Gbf^\ad(\Rbb)^+$, and is an Hermitian symmetric domain.

Take $\Gamma\subset\Gbf(\Rbb)^+$ a congruence subgroup, we have the (connected) Shimura variety $M=\Gamma\bsh X^+$, where $\Gamma$ acts on $X^+$ through its image in $\Gbf^\ad(\Rbb)^+$. The theorem of Baily-Borel compactification affirms that $M$ is a normal quasi-projective algebraic variety over $\Cbb$, and we always assume that $\Gamma$ is torsion-free, so that $M$ is smooth and its  fundamental group is identified with $\Gamma$.

The canonical projection $\wp=\wp_\Gamma:X^+\ra\Gamma\bsh X^+, x\mapsto \Gamma x$ is called the uniformization map. A Shimura subvariety in $M$ is given as $M'=\wp(X'^+)$ for some connected Shimura subdatum $(\Gbf',X';X'^+)$, namely $(\Gbf',X')$ is a Shimura subdatum of $(\Gbf,X)$ in the sense of \cite{deligne pspm} and $X'^+$ is a connected component of $X'$ contained in $X^+$. It is known that $M'$ is a closed subvariety in $M$, and its fundamental group is isomorphic to $\Gamma\cap\Gbf'^\der(\Rbb)^+$.

\end{definition}

\begin{example}[Siegel modular variety]\label{example siegel modular variety}

Fix $V$ a rational symplectic space of dimension $2g$, we have the Shimura datum $(\GSp_V,\Hcal_V;\Hcal_V^+)$ where $\Hcal_V^+$ is the Siegel upper half space of genus $g$. Usually we assume that $V$ comes from the standard symplectic $\Zbb$-module of discriminant 1, and thus a suitable choice of a congruence subgroup $\Gamma\subset\GSp_V(\Qbb)$ defines Siegel modular varieties $\Acal_g=\Acal_V:=\Gamma\bsh \Hcal_V^+$ parameterizing principally polarized abelian variaties of dimension $g$ with level-$\Gamma$ structures. 

We mainly consider Shimura subvarieties in $\Acal_V$. As is explained in \cite{chen lu zuo compositio}, a Shimura subdatum $(\Gbf,X;X^+)$ defining a Shimura subvariety $M\subset\Acal_V$ gives rise to a rational symplectic representation $\Gbf\ra\GSp_{V}$ satisfying Satake's condition (H2) in the sense of \cite{satake rational}. We always assume that the level structure $\Gamma$ is suitably chosen so that the inclusion $M\mono \Acal_g$ extends to their smooth toroidal compactifications $\Mbar\mono\Acalbar_V$, which joins to $M$ resp. to $\Acal_V$ finitely many boundary divisors.

\end{example}

We write $\Tcal_g^\circ$ for the schematic image of $$j:\Mcal_g\ra\Acal_g=\Acal_V,\ [C]\mapsto[\Jac(C)]$$ called the open Torelli locus, and $\Tcal_g$ for its closure, called the Torelli locus. The slope inequality of Xiao is mainly applied to the following situation:

\begin{proposition}[Higgs bundles for surface fibration]\label{proposition higgs bundles for surface fibration} Let $C$ be a closed curve in $\Acalbar_g$ contained generically in the Torelli locus, namely  $C^\circ:=C\cap\Tcal_g^\circ$ is open in $C$. Let $B^\circ$ be the normalization of the preimage of $C^\circ$ in $\Mcal_g$, giving rise to a family of curve $f^\circ:S^\circ\ra B^\circ$ which is compactified into a surface fibration $f:S\ra B$ with semi-stable fibers of genus $g$. Write $i:B\ra C$ for the induced morphism from $B$ into $C$, and $\Vcal_C^{1,0}$ for the $(1,0)$-part of the logarithmic Higgs bundle $\Vcal_C$ on $C$ deduced from the variation of Hodge structure on $\Acal_g$ defined by the moduli problem. Then holds the isomorphism $$f_*\omega_{S/B}\isom i^*\Vcal^{1,0}_C$$ where $\omega_{S/B}$ is the relative dualizing sheaf for $f$. In particular we have the decomposition $f_*\omega_{S/B}=\Fcal_B\oplus\Acal_B$, where $\Fcal_B$ is the semi-stable subbundle of slope 0 corresponding to the Higgs subbundle in $\Vcal_C$ given by the maximal unitary subrepresentation in of the $\Cbb$-linear representation of $\pi_1(C)$ associated to $\Vcal_C$ using Simpson's correspondence.

\end{proposition}

In the rest of this section we derive a naive  bound on the flat part in the canonical Higgs bundle associated to a Shimura variety of dimension $>0$ contained in $\Tcal_g$ generically using Xiao's inequality. 

\begin{proposition}[naive bound on the flat part]\label{proposition naive bound on the flat part} Let $M\subset\Acal_g$ be a Shimura subvariety defined by a subdatum $(\Gbf,X;X^+)$, such that the derived part of the $\Qbb$-group  $\Gbf$ admits an isomorphism $\Gbf^\der=\Res_{F/\Qbb}\Hbf$ for some totally real field $F$ and some semi-simple $F$-group $\Hbf$, and the representation $\Gbf^\der\mono\Sp_V$ decomposes into  $$V=V_0\oplus\Res_{F/\Qbb}W$$ with $V_0$ a trivial subrepresentation and $\Hbf$ acting on $W$ preserving some symplectic $F$-form. Assume further that $F$ is of degree $d$ over $\Qbb$, such that: \begin{itemize}
\item along $r(>0)$ real embeddings $\sigma_1,\cdots,\sigma_r:F\mono\Rbb$, the Lie group $\Hbf(\Rbb,\sigma_i)$ is compact;

\item along the other $d-r(>0)$ real embeddings $\sigma_{r+1},\cdots,\sigma_d:F\mono\Rbb$, the Lie group $\Hbf(\Rbb,\sigma_j)$ is non-compact.
\end{itemize}Here $\Hbf(\Rbb,\sigma_i)$ is the evaluation of $\Hbf$ at $\sigma_i:F\mono\Rbb$.

If $M$ is contained generically in the Torelli locus, then holds the inequality $$\frac{r}{d}\leq \frac{5}{6}+{\frac{1}{6g}}.$$
\end{proposition}

\begin{proof}

The $\Qbb$-linear representation $V=V_0\oplus \Res_{F/\Qbb}W$ decomposes into
 $$V\otimes_\Qbb\Rbb=V_0\otimes_\Qbb\Rbb\oplus\Big(\bigoplus_{i=1}^{d}W\otimes_{F,\sigma_i}\Rbb\Big)$$
after the base change $\Qbb\mono\Rbb$, with $\Gbf^\der(\Rbb)$ acting on $V_0\otimes_\Qbb\Rbb$ trivially,
and $\Hbf(\Rbb,\sigma_i)$ acting on $V\otimes_\Qbb\Rbb$ through the summand $W\otimes_{F,\sigma_i}\Rbb$.
Let $$\Vcal=\Vcal_M=\Vcal_0\oplus\bigoplus_{i=1}^{d}\Wcal_i$$ be the Higgs bundle decomposition on $\Mbar$
according to the decomposition of $$V\otimes_\Qbb\Cbb=(V\otimes_\Qbb\Rbb)\otimes_\Rbb\Cbb,$$
where $\Vcal_0$ corresponds to $V_0$ and $\Wcal_i$ corresponds to $W\otimes_{F,\sigma}\Rbb\otimes_\Rbb\Cbb$ respectively.
Then the Higgs subbundles $\Vcal_0$ and $\Wcal_i$ for $i=1,\cdots,r$ are flat as the $\Gbf^\der(\Rbb)$-action factors through compact Lie groups. In particular, the portion of of flat part in $\Vcal$ is at least
$$\frac{\rank \Fcal^{1,0}}{\rank \Vcal^{1,0}}=\frac{\rank\Fcal}{\rank\Vcal}=\frac{\rank\Vcal_0+\sum\limits_{i=1}^{r}\rank\Wcal_i}{\rank\Vcal}\geq \frac{\sum\limits_{i=1}^{r}\rank\Wcal_i}{\sum\limits_{i=1}^{d}\rank \Wcal_i}=\frac{r}{d}$$
since $\rank\Wcal_i=2\dim_FW$ for all $i$.

Assume that $M$ is contained generically in $\Tcal_g$. Then the linear system of the ample line bundle of top degree automorphic forms on $M$ produces a curve $C$ in $\Mbar$ such that $C\cap M$ is open and dense in $C$, and up to Hecke translation we may assume further that $C$ is contained generically in $\Tcal_g$. Restricting $\Vcal_M$ to $C$ gives us the Higgs bundle $\Vcal_C$, which  flat part $\Fcal_C$ contains the pull-back of $\Fcal$ to $C$, hence $$\frac{\rank \Fcal_C}{\rank\Vcal_C}\geq\frac{\rank \Fcal}{\rank \Vcal}=\frac{r}{d}$$ It remains to notice that the map $i:B\ra C$ is finite, and we have $$\rank\Fcal_B=\rank\Fcal_C$$ for the flat part $\Fcal_B$ in $f_*\omega_{S/B}$, hence Xiao's inequality \autoref{xiao's inequality} would fail for $B$ as long as $\frac{r}{d}>\frac{5}{6}+\frac{1}{6g}$, which is absurd. \end{proof}

\begin{corollary}\label{corollary naive bound}Let $M\subset\Acal_V=\Acal_g$ be a Shimura subvariety defined by $(\Gbf,X;X^+)$ with $\Gbf^\der\isom\Res_{F/\Qbb}\Hbf$ for some semi-simple group $\Hbf$ over a totally real field $F$, which is compact along $r$ real embeddings of $F$ and non-compact along the remaining $d-r$ embeddings, $d$ being the degree $[F:\Qbb]$. If $\frac{r}{d}>{\frac{5}{6}}+\frac{1}{6g}$, then $M$ only meets $\Tcal_g^\circ$ at finitely many points.
\end{corollary}

The proof is similar and immediate: otherwise one finds a curve in $\ol M$ contradicting Xiao's inequality. % otherwise.

\section{Shimura subvarieties of $\SU(n,1)$-type} We first briefly recall the set-up for Shimura subvarieties of $\SU(n,1)$-type, which is slightly more general than the one used in \cite{chen lu zuo compositio}, as we no longer require the existence of an Hermitian form over a CM field.

\begin{definition}[Shimura subvarieties of $\SU(n,1)$-type]\label{definition shimura subvarieties of su(n,1)-type} For the fixed symplectic $\Qbb$-vector space $V$ of dimension $2g$, a Shimura subvariety of $\Acal_V$ is said to be of $\SU(n,1)$-type if it is defined by a Shimura subdatum $(\Gbf,X;X^+)\subset(\GSp_V,\Hcal_V;\Hcal_V^+)$ such that $\Gbf^\der\isom\Res_{F/\Qbb}\Hbf$ with \begin{itemize}
	\item $F$ is a totally real number field of degree $d$ and $\Hbf$ is a simple $F$-group;
	
	\item among the real embeddings $\sigma_1,\cdots,\sigma_d$, we have $\Hbf(\Rbb,\sigma_1)\isom\SU(n,1)$, and $\Hbf(\Rbb,\sigma_i)\isom\SU(n+1)$ for $i=2,\cdots,d$.
\end{itemize} We call $F$ the definition field of the datum.
\end{definition}
Note  that when $n\geq2$, $X=X^+$ is connected. The inclusion
$$(\Gbf,X;X^+)\mono(\GSp_V,\Hcal_V;\Hcal_V^+)$$
gives rise to the representation of $\Gbf$ on $V$,
and from \cite{satake rational} we know that the restriction of this representation to $\Gbf^\der=\Res_{F/\Qbb}\Hbf$
is decomposed into $V=V_0\oplus V'$ where
\begin{itemize}
	\item $\Gbf^\der$ acts on $V_0$ trivially;
	
	\item the representation $\Gbf^\der\ra\GL_{V'}$ is the scalar restriction of an $F$-linear representation $\Hbf\ra\GL_{W,F}$ where $W$ is an $F$-vector space carrying a symplectic $F$-form preserved by $\Hbf$.
\end{itemize}

Let $M$ be the Shimura subvariety in $\Acal_V$ of $\SU(n,1)$-type in the sense above, defined by $(\Gbf,X;X)$, isomorphic to $\Gamma\bsh X$ for a torsion-free congruence subgroup $\Gamma\subset\Gbf^\der(\Qbb)^+$, and we write $\Vcal$ for the Higgs bundle  of the $\Cbb$-PVHS on $M$ given by the moduli problem. Similar to the situation in \cite{chen lu zuo compositio} which we have also used in Section 2 for the naive bound, there exists a decomposition of Higgs bundles $$\Vcal=\Vcal_0\oplus\sum_{i=1}^{d}\Wcal_i$$ with $\Vcal_0$ a trivial Higgs bundle corresponding to $V_0$, and $\Wcal_i$ the Higgs bundle associated to $W\otimes_{F,\sigma_i}\Rbb$, in which the only non-flat part is given by $\Wcal_1$. 

According to the Satake classification over $\Rbb$ cf.\cite{satake real},  along $\sigma_1:F\mono\Rbb$ the base change admits a decomposition $$\sigma_1^*W=W\otimes_{F,\sigma_1}\Rbb\isom\bigoplus_{m=1}^{n-1}\Lambda_m^{\oplus r_m}$$ where $\Lambda_m=\wedge^m_\Cbb\Std$ is the $m$-th exterior power of the standard representation of $\SU(n,1)$ on $\Cbb^{n+1}$. Note that the action of $\SU(n,1)$ on $\Lambda_m$ preserves an Hermitian form of signature $\Big(\binom{n}{m},\binom{n}{m-1}\Big)$, and $\dim_\Rbb\Lambda_m=2\binom{n+1}{m}$.

In the sequel we assume  for simplicity that $V$ is primary of type $\Lambda_m$ in the sense of \cite[Subsection 5.2]{chen lu zuo compositio}, namely $V_0=0$ and $\sigma_1^*W\isom(\wedge^m\Std)^{\oplus r}$ for some multiplicity $r\geq 1$.%, with $\Lambda_m$ the symplectic representation of $\SU(n,1)$ on the $m$-th wedge product $\wedge^m(\Std)$, $m=1,\cdots,n$.

We need the following fact to compute the Harder-Narasimhan filtration on certain curves, cf. \cite[7.12]{moeller viehweg zuo}:

\begin{lemma}[decomposition of the canonical Higgs bundle]\label{lemma decomposition of the canonical higgs bundle} Let $\Ecal$ be the Higgs bundle on $M$ associated to the $\Cbb$-representation $$\Gamma\mono\Gbf^\der(\Rbb)\ra\Hbf(\Rbb,\sigma_1)\overset{\Std}{\ra}\GL_\Cbb(\Cbb^{n+1}\otimes_\Rbb\Cbb)$$ (through the unique non-compact factor of $\Gbf^\der(\Rbb)$). Then the Hodge decomposition of $\Ecal$ is of the form $\Ecal=\Ecal^{1,0}\oplus\Ecal^{0,1}$ with \begin{itemize}
		\item $\Ecal^{1,0}\isom\Omega^1_M\otimes\Lcal^\vee\oplus\Lcal$;
		\item $\Ecal^{0,1}\isom(\Ecal^{1,0})^\vee=\Tan_M\otimes\Lcal\oplus\Lcal^\vee$.
	\end{itemize} Here $\Lcal$ is a line bundle on $M$ such that $\Lcal^{\otimes(n+1)}=\omega_M(=\Omega_M^n)$, and the summands in $\Ecal^{1,0}$ and $\Ecal^{0,1}$ are stable subbundles.
	
	In particular, writing $\Scal=\Omega^1_M\otimes\Lcal^\vee$ which gives $\Ecal^{1,0}=\Scal\oplus\Lcal$, we have $$c_1(\Lcal)=\frac{1}{n+1}c_1(\omega_M)=c_1\big(\Omega^1_M\otimes\Lcal^\vee\big),\qquad c_1(\Ecal^{1,0})=2c_1(\Lcal)=\frac{2}{n+1}c_1(\omega_M).$$
	For the $m$-th exterior power $$\wedge^m\Ecal^{1,0}=\wedge^m\Scal\oplus\wedge^{m-1}\Scal\otimes\Lcal,$$
	we have
	$$\begin{aligned}
	c_1(\wedge^m\Scal)&\,=\binom{n-1}{m-1}\cdot \frac{1}{n+1}c_1(\omega_M),\\
	c_1(\wedge^{m-1}\Scal\otimes\Lcal)&\,=\Bigg(2\binom{n}{m-1}-\binom{n-1}{m-1}\Bigg)\cdot \frac{1}{n+1}c_1(\omega_M),\\
	c_1(\wedge^m\Ecal^{1,0})&\,=\binom{n}{m-1}\cdot \frac{2}{n+1}c_1(\omega_M).
	\end{aligned}$$
	
	Similar properties hold for the canonical logarithmic extensions of these structures to $\Mbar$ (the smooth toroidal compactification of $M$).%, equal to the closure of $M$ in a suitable smooth toroidal compactification of $\Acal_g$).
	
\end{lemma}

\begin{proof}
The decomposition for $\Ecal$ and the first Chern classes for $\Scal$ and $\Lcal$ are given in \cite{moeller viehweg zuo}. The claim for higher exterior powers follows from the splitting principle of Chern classes
(using the facts that $c_1(\wedge^m\Mcal)=\binom{\rank(\Mcal)-1}{m-1}\cdot c_1(\Mcal)$ and $c_1(\Mcal'\otimes\Mcal'')=\rank(\Mcal')c_1(\Mcal'')+\rank(\Mcal'')c_1(\Mcal')$\big).
\end{proof}

\begin{theorem}[exclusion of Shimura subvarieties of $\SU(n,1)$-type]\label{theorem exclusion of shimura subvarieties of su(n,1)-type}Let $M\subset\Acal_V=\Acal_g$ be a Shimura subvariety of $\SU(n,1)$-type defined by a Shimura subdatum $(\Gbf,X;X)$ as in \autoref{definition shimura subvarieties of su(n,1)-type}, such that the rational symplectic representation $\Gbf^\der\mono\Sp_W$ is primary of type $\Lambda_m$ in the sense of \cite{chen lu zuo compositio}.%, namely it is restricted from some homomorphism of $F$-groups $\Hbf\mono\Sp_V$ for some symplectic $F$-form $V$ such that $$\sigma_1^*V:=V\otimes_{F,\sigma_1}\Rbb\isom(\wedge^m\Std)^{\oplus r}$$ as real symplectic representations of $\Hbf(\Rbb,\sigma_1)\isom\SU(n,1)$. 
	Then $M$ is NOT contained generically in the Torelli locus in $\Acal_W$ as long as the following inequality holds: $$\frac{n+m-1}{n}\left(\frac{n+1}{m}\cdot d -\frac{2}{r\binom{n}{m-1}}\right)>12.$$
	
\end{theorem}

\begin{proof}

Assume on the contrary that $M$ is contained generically in $\Tcal_g$. Then the same construction used in \autoref{proposition naive bound on the flat part} produces a curve $C$ in the $n$-dimensional subvariety $\Mcalbar(\subset\Acalbar_V)$ using successive $n-1$ hyperplane sections of a fixed very ample power $\omega_M^N$ of the automorphic line bundle on $M$, such that $C^\circ=C\cap\Tcal_g^\circ$ is open and dense in $C$, which gives rise to a surface fibration $f:S\ra B$ and a finite morphism $i:B\ra C$. 

The logarithmic Higgs bundle $\Vcal^{1,0}$ on $\Mbar$ defined by the moduli problem $M\mono\Acal_V$ already admits the following filtration $$0=\Vcal_0\subsetneq\Vcal_1\subsetneq\Vcal_2\subsetneq\Vcal_3=\Vcal^{1,0}$$ where: \begin{itemize}
\item $\Vcal_1$ is the summand $(\wedge^{m-1}\Scal\otimes\Lcal)^{\oplus r}$ in $\Wcal_1$;

\item $\Vcal_2=\Wcal_1\isom\Vcal_1\oplus(\wedge^m\Scal)^{\oplus r}$;

\item $\Vcal_3=\Vcal^{1,0}$ only differs from $\Wcal_1$ by a direct summand $\Fcal$ flat of degree zero.
\end{itemize}

The graded quotients of this filtration are already semi-stable. Computing their degrees over $C$ we find that the maximal slope in $\Vcal^{1,0}_C$ is realized on $\Vcal_1$
$$\mu(\Vcal_1)=\mu(\wedge^{m-1}\Scal\otimes\Lcal)=\frac{2\binom{n}{m-1}-\binom{n-1}{m-1}}{\binom{n}{m-1}}\cdot d_C$$
with the constant $d_C=\frac{N^{n-1}}{n+1}c_1(\omega_\Mbar)^n$
as we have used the fixed very ample line bundle $\omega^N_M$ on $M$ which extends to $\omega_\Mbar^N$ on $\Mbar$.   
Note also that
$$\begin{aligned}
\deg(\Vcal^{1,0}_C)&\,=r\deg(\wedge^m\Ecal^{1,0})=2r\binom{n}{m-1}\cdot d_C,\\
g&\,=\frac{1}{2}\dim_\Qbb V=rd\binom{n+1}{m},
\end{aligned}$$
and that the maximal slope part is $\Vcal_1\isom(\wedge^{m-1}\Scal\otimes\Lcal)^{\oplus r}$, which is of rank $r\binom{n}{m-1}$. %$\rank\Fcal_C=\frac{d-1}{d}g$.

Passing from $C$ to $B$ only changes $d_C$ (resp. $\mu(\Vcal_1)$) to $d_B$ (resp. $\mu(i^*\Vcal_1)$)
  by a positive multiple $\deg i$, and Xiao's inequality gives us
  $$12\cdot 2r\binom{n}{m-1}\cdot d_B\geq\Bigg(2rd\binom{n+1}{m}-2+r\binom{n}{m-1}\Bigg)\cdot
  \Bigg(2-\frac{\binom{n-1}{m-1}}{\binom{n}{m-1}}\Bigg)\cdot d_B,$$
   which is $$\frac{n+m-1}{n}\left(\frac{n+1}{m}d+\frac{1}{2}-\frac{1}{r\binom{n}{m-1}}\right)\leq 12$$
   and hence the generic exclusion when the inequality above fails.\end{proof}
% $$\frac{n+m-1}{n}(\frac{n+1}{m}(d+1)-\frac{2}{r\binom{n}{m-1}})\leq 12$$ and hence the generic exclusion when the inequality above fails. \end{proof}

For example, when $m=1$, we obtain the generic exclusion for such Shimura subvarieties satisfying $$(n+1)d>\frac{23}{2}+\frac{1}{r}$$ which is \begin{itemize} 
\item $(n+1)d\geq 13$ when $r=1,2$; this is, for example, the case when the symplectic representation is directly obtained as the imaginary part of an Hermitian form of signature $(n,1)$ over some CM quadratic extension $E/F$, and $g=d(n+1)$;

\item $(n+1)d\geq 12$ when $r\geq 3$; this is the same as the inequality $d(n+1)\geq 12$ given in \cite{chen lu zuo compositio} under more restrictive assumptions.

\end{itemize}
% $$(d+1)(n+1)>12+\frac{2}{r}$$ which is \begin{itemize}
	%	\item $(d+1)(n+1)\geq 15$ when $r=1$: this is, for example, the case when the symplectic representation is directly obtained as the imaginary part of an Hermitian form of signature $(n,1)$ over some CM quadratic extension $E$, and $g=d(n+1)\geq 15-n-1$.
		
	%	\item $(d+1)(n+1)\geq 13$ when $r\geq 2$: this automatically becomes $(d+1)(n+1)\geq 14$ because 13 is prime while $d+1$ and $n+1$ are both at least 2, and this is slightly finer than the inequality $d(n+1)\geq 12$ given in \cite{chen lu zuo compositio} Corollary 1.2.5 with  $m=1$. % for $n$ large.
	%\end{itemize}

\section{Shimura varieties of orthogonal type}

In this section we consider a class of Shimura varieties of orthogonal type
containing Shimura subvarieties of $\SU(n,1)$-type. Recall the following:

\begin{definition}[Shimura subvarieties of orthogonal type] A Shimura subvariety of orthogonal type, or more precisely, of $\Spin(N,2)$-type, in $\Acal_V$ is defined by a subdatum $(\Gbf,X;X^+)$ such that $\Gbf^\der(\Rbb)\isom\Spin(N,2)\times\Spin(N+2)^{d-1}$ for some $d$. Here $\Spin(a,b)$ is the spin group of the standard quadratic space $\Rbb^{a+b}$ of signature $(a,b)$.

The Hermitian symmetric domain $X$ above is the one associated to $\Spin(N,2)$, and we often make use of the following two equivalent descriptions of $X$: \begin{itemize}
\item[(1)] $X$ is the open subset of two-dimensional negative definite $\Rbb$-subspaces in $\Rbb^{N+2}$ in the Grassmannian $\Grbf(2,\Rbb^{N+2})$;

\item[(2)] $X$ is the open subset of negative definite isotropic $\Cbb$-lines in $\Pbb(\Cbb^{N+2})$, namely those $\Cbb v$ such that $b(v,\vbar)=0$ and $b(v,\vbar)<0$ for $b$ the quadratic form of signature on $\Rbb^{N+2}$ extended to $\Cbb^{N+2}$; in particular $X$ is an open subset of the quadric defined by $b$ in $\Pbb(\Cbb^{N+2})$.\end{itemize}

The equivalence between (1) and (2) is well-known: starting with $u,u'\in \Rbb^{N+2}$ which are negative definite and orthogonal to each other, we can choose a suitable $J\in\Cbb^\times$ purely imaginary so that $v=u+Ju'$ defines a line in (2); conversely, given a line $\Cbb v$ in (2), one produces a pair of negative definite vectors $(u,u')$ orthogonal to each other giving rise to a negative definite $\Rbb$-plane in (1).

%In fact,  for  $b$ the quadratic form on $\Rbb^{N+2}$ of signature $(N,2)$ which extends to $\Cbb^{N+2}$, $X$ can be identified with the set of isotropic $\Cbb$-lines in $\Pbb(\Cbb^{N+2})$ which are negative definite, i.e. $\Cbb v$ such that $b(v,v)=0$ and $b(v,\vbar)<0$. 

\end{definition}

The two descriptions above also give us a natural equivariant embedding of Hermitian symmetric domains: if $(V,q)$ is a quadratic space of signature $(N,2)$ over $\Rbb$, and $U\subset V$ is a positive definite subspace of dimension $N-2n(\geq 0)$, such that the restriction of $q$ to the orthogonal complement of $U$ in $V$ is the real part of some Hermitian space $(W,h)$ of signature $(n,1)$, then we have a natural inclusion of semi-simple Lie groups $\SU(W,h)\mono\SO(V,q)$ inducing an equivariant holomorphic embedding of the corresponding Hermitian symmetric domains $X(W,h)\mono X(V,q)$, sending a negative definite $\Cbb$-line in $X(W,h)$ to the associated negative definite $\Rbb$-plane in $X(V,q)$. Note that this inclusion actually factors through $$\SU(W,h)\mono\SO(W,b)\mono\SO(V,q)$$ and $$X(W,h)\mono X(W,b)\mono X(V,q),$$
where we write $b$ for the real part of $h$, equal to the restriction of $q$ to $W$ as a real subspace of $V$. Also the inclusions $\SU(W,h)$ into special orthogonal groups  lift into homomorphisms into the corresponding spin groups, because special unitary groups are simply connected.

In this section we are only interested in the Coleman-Oort problem for Shimura varieties of $\Spin(N,2)$-type containing Shimura varieties of $\SU(n,1)$-type. We start with the case $N=2n$ before entering the general case where $N> 2n$.

Let $(W,h)$ be an Hermitian space over a CM field $E$, and write $F$ for the totally real part of $E$. Assume that the signature of $h$ is: \begin{itemize}
\item $(n,1)$ along one real embedding $\sigma:F\mono\Rbb$;

\item definite along the other embeddings $\sigma_2,\cdots,\sigma_d:F\mono\Rbb$
\end{itemize} with $d=[F:\Qbb]$. We also have the real part of $(W,h)$, namely $(U,q)$ with $U=\Res_{E/F}W$ and $q=\tr_{E/F}h$. Write $\Hbf=\SU(W,h)$ for the special unitary $F$-group associated to $(W,h)$, contained in $\Hbf'=\Spin(U,q)$ the spin $F$-group of $(U,q)$. Let $M'\subset\Acal_V=\Acal_g$ be a Shimura subvariety of $\Spin(2n,2)$-type, i.e. defined by a Shimura subdatum $(\Gbf',X';X'^+)$ with $\Gbf'^\der=\Res_{F/\Qbb}\Hbf'$, which contains a Shimura subvariety $M$ of $\SU(n,1)$-type defined by a Shimura subdatum $(\Gbf,X;X^+)$ with $\Gbf^\der=\Res_{F/\Qbb}\Hbf$. 

We have used spinor representations  in \cite{chen lu zuo compositio}: \begin{itemize}

\item the spinor representations of the Lie group $\Spin(2n,2)$  are $P_+$ and $P_-$: choose any splitting of $\Cbb^{2n+2}$ into $T\oplus T^\vee$ such that the quadratic form is equivalent to the pairing $$((u,u^\vee),(v,v^\vee))\mapsto u^\vee(v)+v^\vee(u)$$ we have $P_+=\wedge^+T:=\wedge^\even_\Cbb T$ and $P_-=\wedge^-T:=\wedge^\odd_\Cbb T\isom P_+^\vee$, both of dimension $2^{n}$ over $\Cbb$.

\end{itemize}

We need to find suitable splitting for $(U,q)$, at least over $\Cbb$, using $(W,h)$. Consider the following lemma for the quadratic extension $\Cbb/\Rbb$, in which we temporarily use $(U,q)$ and $(W,h)$ to denote the quadratic and Hermitian forms involved:

\begin{lemma}
Let $h:W\times W\ra\Cbb$ be an Hermitian form over $\Cbb$, with real part $q:U\times U\ra\Rbb$, and $q_\Cbb:U_\Cbb\times U_\Cbb\ra\Cbb$. Then $U_\Cbb$ admits a natural splitting by $W\mono U_\Cbb=W\otimes_\Rbb\Cbb$, and $U_\Cbb\isom W\oplus\Wbar$ with respect to the real structure on $U_\Cbb$ given by $U$.
\end{lemma}

\begin{proof}
We may diagonalize $(W,h)$ into direct sums of one dimensional $\Cbb$-spaces, and reduce to the case when $\dim_\Cbb W=1$: $W=\Cbb w$ for some basis $w$, and $h(zw,z'w)=c\zbar z'$ for $z,z'\in\Cbb$ and $c\in\Rbb^\times$. Thus $U$ admits a basis $(w,\ibf w)$ for a fixed choice $\ibf=\sqrt{-1}$, and $q=\tr_{\Cbb/\Rbb}h$ sends $(aw+b\ibf w,a'w+b'\ibf w)$ to $c(aa'+bb')$. It suffices to choose the splitting to be $$W\mono U_\Cbb=W\otimes_\Rbb\Cbb,\ w\mapsto w\otimes(1+\ibf)
$$ in which case $\Wbar$ is simply $W\otimes(1-\ibf)$.\end{proof}

Return to the general setting over a CM field $E/F$. Evaluate the inclusion $\Hbf=\SU(W,h)\mono\Hbf'=\Spin(U,q)$ at $\sigma=\sigma_1:F\mono\Rbb$, we get $\SU(n,1)\mono\Spin(2n,2)$. The $\Cbb$-representations $P_\pm\isom \wedge^\pm(W_\Cbb)$ of $\Hbf'(\Rbb,\sigma)$ restricts to sums of wedge product $\Cbb$-representations of $\Hbf(\Rbb,\sigma)$:

\begin{lemma}

Let $M\mono M'$ be the inclusion of a Shimura variety of $\SU(n,1)$-type into a Shimura variety of $\Spin(2n,2)$-type associated to $\Hbf\mono\Hbf'$ as above. Consider the Higgs bundles $\Pcal_\pm$ on $M'$ associated to the spinor representation $P_\pm$ of $\Hbf'(\Rbb,\sigma)$. Then the restriction of $\Pcal_\pm$ to $M$ decomposes into $$\Pcal_+=\bigoplus_{m\ \even}\Ecal_m,\ \Pcal_-=\bigoplus_{m\ \odd}\Ecal_m$$ with $\Ecal_m$ the Higgs bundle associated to the $m$-th exterior power $\Lambda_m$ of the the representation $$\Gbf^\der(\Rbb)\ra\Hbf(\Rbb,\sigma_1)\overset{\Std}{\lra}\GL_{n+1}(\Cbb).$$
%we have a decomposition of the corresponding Higgs bundles associated to the $\Cbb$-representations of $\pi_1(M)$ on $P_\pm$ through $\Gbf^\der(\Rbb)\ra \Hbf(\Rbb,\sigma)\mono\Hbf'(\Rbb,\sigma)$: $$\Pcal_+=\bigoplus_{2|m}\Ecal_m,\ \Pcal_-=\bigoplus_{2\not|i}\Ecal_m$$ where $\Pcal_\pm$ is the Higgs bundle on $M$ corresponding to $P_\pm$, and $\Ecal_m$ corresponds to the $m$-th exterior power $\Lambda_m$ of the standard representation $\Std=\Cbb^{n+1}$ for $\Hbf(\Rbb,\sigma)=\SU(n,1)$.

In particular, $\Ecal_1=\Ecal$ and $\Ecal^{1,0}=\Scal\oplus\Lcal$, and $\Ecal_m^{1,0}$ decomposes into the direct sum of $\wedge^m\Scal$ and $\wedge^{m-1}\Scal\otimes\Lcal$, whose ranks and Chern classes are given as in \autoref{lemma decomposition of the canonical higgs bundle}. Take $C$ a generic curve in $\Mbar$ produced from the linear system of $\omega_M^N$ a fixed very ample power of $\omega_M$,  the slopes of these summands on $C$ are as follows: $$\mu(\wedge^m\Scal_C)=\frac{m}{n}d_C,\quad \mu(\wedge^{m-1}\Scal_C\otimes\Lcal_C)=\frac{2n-m}{n}d_C$$ for $m=1,\cdots,n$, with  $d_C=\frac{N^{n-1}}{n+1}c_1(\omega_\Mbar)^n$. We also have $c_1(\Ecal_0^{1,0})=0$ and $c_1(\Ecal_{n+1}^{1,0})=\frac{2}{n+1}c_1(\omega_M)$, and thus $\mu(\Ecal_{0,C}^{1,0})=0$ and $\mu(\Ecal_{n+1,C}^{1,0})=2d_C$.

\end{lemma}The proof is immediate after \autoref{lemma decomposition of the canonical higgs bundle}.

The following proposition takes care of the inclusion $\SU(n,1)\mono\Spin(2n,2)$ in a natural arithmetic setting. Note that the case (1) is singled out for later use in \autoref{theorem generic exclusion of Shimura varieties of orthogonal types}.

\begin{proposition}[restriction on the real part]\label{proposition restriction on the real part} Let $M'\subset\Acal_V=\Acal_g$ be a Shimura subvariety of $\Spin(2n,2)$-type, containing a Shimura subvariety $M\subset M'$ of $\SU(n,1)$-type in the sense above. Assume that the symplectic representation $\Gbf'^\der\mono\Sp_V$ defining $M'\mono\Acal_V$ admits no trivial subrepresentations. If $M'$ is contained generically in $\Tcal_g$, then the following hold:

(1) if $P_+$ and $P_-$ appear with equal multiplicity $r$, then $d\leq 3+2^{-n-1}(\frac{1}{r}-\frac{1}{2})$;

(2) if $M'\mono\Acal_V$ is primary of type $P_-$ of multiplicity $r$, then along the parity of $n$ we have: \begin{itemize}

\item $d\leq 3+2^{-n}(\frac{1}{r}-\frac{1}{2})$ when $n$ is even;

\item $d\leq \frac{6n}{2n-1}+2^{-n}(\frac{1}{r}-\frac{1}{2})$ when $n$ is odd;

\end{itemize}

(3) if $M'\mono\Acal_V$ is primary of type $P_+$ of multiplicity $r$, then along the parity of $n$ we have: \begin{itemize}
\item $d\leq \frac{3n}{n-1}+2^{-n}(\frac{1}{r}-\frac{n}{2})$ when $n$ is even;

\item $d\leq 3+2^{-n}(\frac{1}{r}-\frac{1}{2})$  when $n$ is odd.

\end{itemize}

\end{proposition}

\begin{proof}

We have $\Gbf^\der=\Res_{F/\Qbb}\Hbf$ and $\Gbf'^\der=\Res_{F/\Qbb}\Hbf'$, and the inclusion $\Gbf'\mono\Sp_V$ is restricted from $\Hbf'\mono\Sp_W$ where $W$ is an $F$-symplectic space, so that $V\isom\Res_{F/\Qbb}W$ and the canonical logarithmic Higgs bundle $\Vcal=\Vcal_{M'}$ admits a decomposition $\Vcal=\oplus_{i=1,\cdots,d}\Wcal_i$ with $\Wcal_i$ corresponds to the action of $\Gbf'^\der(\Rbb)$ on $W\otimes_{F,\sigma_i}\Rbb$ through $\Hbf'(\Rbb,\sigma_i)$. In particular, the $(1,0)$-parts $\Wcal_2^{1,0},\cdots,\Wcal_d^{1,0}$ are already flat.

Similar to \autoref{theorem exclusion of shimura subvarieties of su(n,1)-type}, we may assume that a generic curve $C$ is chosen in $\Mbar'$ such that $C$ is contained generically in $\Tcal_g$ and lifts to a semi-stable surface fibration $f:S\ra B$ together with a finite morphism $i:B\ra S$ satisfying \autoref{proposition higgs bundles for surface fibration}.

(1) In this case we have $$W\otimes_{F,\sigma_1}\Rbb\isom(P_+\oplus P_-)^{\oplus r}\isom(\wedge^\cdot(\Cbb^{n+1}))^{\oplus r}$$ for some multiplicity $r$, and thus $\Wcal_1^{1,0}\isom(\oplus_{m=0}^{n+1}\Ecal_m^{1,0})^{\oplus r}$. Therefore the maximal slope in $\Vcal_C^{1,0}$ the restriction to $C$ is given by $$\mu_\maxx=\mu((\Ecal_{n+1,C}^{1,0})^{\oplus r})=\mu(\Ecal_{n+1,C}^{1,0})=2d_C$$ with $(\Ecal^{1,0}_{n+1,C})^{\oplus r}$ of rank $r$, while $\deg\Vcal_C^{1,0}=\deg\Wcal_{1,C}^{1,0}=r\sum_m\deg(\Ecal_{m,C}^{1,0})=r\cdot 2^{n+1}d_C$.  and the flat part is $$(\wedge^0\Ecal_1^{1,0})^{\oplus r}\oplus\Wcal_2^{1,0}\oplus\cdots\oplus\Wcal_d^{1,0}.$$ Here we have used the same constant $d_C$ as in \autoref{theorem exclusion of shimura subvarieties of su(n,1)-type}. Passing to $B$ using $i:B\ra C$, the slopes and degrees only differ by a common multiple $\deg(i)$, namely one replaces the constant $d_C$ by $d_B=\deg(i)d_C$. 

When $M'$ is contained generically in $\Tcal_g$, we may choose $C$ in $\Mbar\subset\Mbar'$ such that $C$ is contained generically in $\Tcal_g$, so that Xiao's inequality for $f:S\ra B$ gives $$12\cdot 2^{n+1}r\geq(2\cdot 2^{n+1}rd-2+r)2$$ which is $d\leq 3+2^{-n-1}(\frac{1}{r}-\frac{1}{2})$.

% $$12\cdot r\cdot 2^{n+1}\geq 2(4rd\cdot 2^{n+1}-4-2r(1+(d-1)2^{n+1}))$$ which is $d\leq 2+\frac{1}{2^nr}+\frac{1}{2^{n+1}}$, namely $d\leq 2$.

(2) In this case $\Vcal_C^{1,0}=\Wcal_{1,C}^{1,0}\oplus\cdots\Wcal_{d,C}^{1,0}$ is of rank $rd\cdot 2^{n}$, with $\Wcal_{2,C}^{1,0},\cdots,\Wcal_{d,C}^{1,0}$ flat, and $\Wcal^{1,0}_C\isom(\Pcal_-^{1,0})^{\oplus r}$ contains no flat part, and it already contains $(\Ecal_{1,C}^{1,0})^{\oplus r}$ of slope $(2-\frac{1}{n})d_C$, which is:

(2-1) maximal when $\Ecal_{n+1,C}$ does not contribute, namely $n$ is odd, and in this case the maximal slope is realized on $\Lcal^{\oplus r}$ contained in $(\Ecal^{1,0}_{1,C})^{\oplus r}$, of rank $r$;

(2-2) strictly smaller than the maximal slope $2d_C=\mu(\Ecal^{1,0}_{n+1,C})$ when $n$ is even, and the maximal slope is realized on $(\Ecal^{1,0}_{n+1,C})^{\oplus r}$ of rank $r$.

Note that $\Wcal_C^{1,0}\isom(\Pcal_-^{1,0})^{\oplus r}$ is of rank $2^nr$ and degree $2^nrd_C$, hence Xiao's inequality implies:

(2-1) when $n$ is odd: $12\cdot 2^nr\geq(2^{n+1}rd-2+r)(2-\frac{1}{n})$, which is $d\leq\frac{6n}{2n-1}+2^{-n}(\frac{1}{r}-\frac{1}{2})$; % $(2-\frac{1}{n})(d+2^{-n}(\frac{1}{2}-\frac{1}{r}))\leq 6$;

(2-2) when $n$ is even: $12\cdot 2^nr\geq(2^{n+1}rd-2+r)\cdot 2$, which is $d\leq 3+2^{-n}(\frac{1}{r}-\frac{1}{2})$.

%(2-1) $6\geq(d+1-\frac{1}{2^{n-1}r})(2-\frac{1}{n})$, namely $d\leq \frac{4n+1}{2n-1}+\frac{1}{2^{n-1}r}$ when $n$ is odd;

%(2-2) $d\leq 2+\frac{1}{2^{n-1}r}$ when $n$ is even.

(3) In this case $\Wcal_C^{1,0}\isom(\Pcal_+^{1,0})^{\oplus r}$ is of rank $2^nr$. Along the parity of $n$ we have: % and contains a trivial subbundle of rank $r$ given by $(\Ecal_{0,C}^{1,0})^{\oplus r}$, which means the flat part in $\Vcal_C$ is of rank $r(1+(d-1)2^n)$. Hence along the parity of $n$ we have:

(3-1) when $n$ is even: the maximal slope comes from the summand $\Scal_C\otimes\Lcal_C$ in $\Ecal_{2,C}^{1,0}$, which is $(2-\frac{2}{n})d_C$ of rank $rn$, and  Xiao's inequality becomes 
$$12\cdot 2^nr\geq(2^{n+1}rd-2+rn)(1-\frac{1}{n})\cdot 2,$$ namely $d\leq \frac{3n}{n-1}+2^{-n}(\frac{1}{r}-\frac{n}{2})$ (note that $n\geq 2$).
% $$d\leq\frac{2n+1}{n-1}+\frac{1}{2^{n-1}r}+\frac{1}{2^n}$$

(3-2) when $n$ is odd: the maximal slope comes from the summand $\Ecal^{1,0}_{n+1,C}$, which is $2d_C$ of rank $r$, and similar to the case (2-2), Xiao's inequality is $12\cdot 2^nr\geq(2^{n+1}rd-2+r)\cdot 2$, namely $d\leq 3+2^{-n}(\frac{1}{r}-\frac{1}{2})$.% $$d\leq 2+\frac{1}{2^{n-1}r}+\frac{1}{2^n}.$$

\end{proof}

In particular, when $r\geq 3$, the generic exclusion of $M$ holds whenever $d\geq\frac{3n}{n-1}$ for $n$ even and $d\geq\frac{6n}{2n-1}$ for $n$ odd, which are slightly finer that the bound $d\geq 6$ given in \cite{chen lu zuo compositio}.

%In particular, the above discussion implies:\begin{itemize}
%\item when $n=1$: $d\leq 5$ if the symplectic representation is primary of $P_-$-type, and $d\leq 3$ if primary of $P_+$-type;

%\item when $n=2$: $d\leq 2$ if primary of $P_-$-type, $d\leq 5$ if primary of $P_+$-type;

%\item when $n\geq 3$: $d\leq 2$ suffices whatever the primary type is.
%\end{itemize}

We proceed to the remaining case when $N>2n$:

\begin{theorem}[generic exclusion of Shimura varieties of orthogonal types]\label{theorem generic exclusion of Shimura varieties of orthogonal types} Let $M'\subset \Acal_V$ be a Shimura variety of $\Spin(N,2)$-type which contains a Shimura subvariety $M$ of $\SU(n,1)$-type in the sense above, namely $M'$ is associated to some quadratic space $(W,q)$ over some CM field $E/F$ subject to the constraints of signature used in this section, and let $U$ be a positive definite subspace of signature $(N-2n,0)$, whose orthogonal complement is the real part of some Hermitian space $(H,h)$ over $E/F$ giving rise to the Shimura subvariety $M$ mentioned above.

Assume that the symplectic representation defining $M'\mono\Acal_V$ is primary. Then $M'$ is NOT contained generically in $\Tcal_g$ as long as $[F:\Qbb]>3+\frac{1}{m\cdot 2^{\lfloor(N+1)/2\rfloor}}-\frac{1}{2^{n+2}}$, where $m$ is the multiplicity of the spinor representation in $W\otimes_{F,\sigma}\Rbb$ for the group $\Spin(N,2)$.

\end{theorem}

\begin{proof}

Recall the behavior of spinor representations to spin subgroups used in \cite{chen lu zuo compositio}, where we consider the case over $\Cbb$ for simplicity:\begin{itemize}

\item the restriction of the spinor representation $P$ of $\Spin(2k+1)$ to $\Spin(2k)$, using an orthogonal decomposition of the form $\Cbb^{2k+1}=L\oplus\Cbb^{2k}$ for some line $L$, is the direct sum of $P_+$ and $P_-$, the two half-spin representations of $\Spin(2k)$;

\item similarly, the restrictions of $P_+$ and $P_-$ of $\Spin(2k+2)$ to $\Spin(2k+1)$
are both isomorphic to the spinor representation of $\Spin(2k+1)$ as long as the restriction comes from an orthogonal decomposition $\Cbb^{2k+2}=L\oplus\Cbb^{2k+1}$ by some line $L$.
\end{itemize}

In our setting, the condition $N-2n>0$ refines the inclusion $W\supset\Res_{E/F}H$ into $W\supset W'\supsetneq\Res_{E/F}H$ with $W'$ of dimension $3+2n$ over $F$, such that \begin{itemize}
\item $W=U'\oplus W'$ is an orthogonal direct sum decomposition with $U'$ of signature $(N-2n-1,0)$ along $\sigma$;

\item $W'=L\oplus \Res_{E/F}H$ is an orthogonal direct sum decomposition with $L$ a positive definite line.
\end{itemize}

The inclusion of Shimura varieties $M\subset M'$ is thus refined into $M\subset M''\subset M'$ with $M''$ the Shimura variety of $\Spin(2n+1,2)$-type associated to $(W',q|_{W'})$, and we write $(\Gbf'',X'';X''^+)$ for the corresponding Shimura subdatum of $(\GSp_V,\Hcal_V;\Hcal_V^+)$. Since $M'\mono\Acal_V$ is defined by some symplectic representation $\Gbf'^\der\mono\Sp_V$ primary of spinor type, its restriction to $\Gbf''^\der\mono\Sp_V$ is primary of spinor type, namely it is the scalar restriction from $\rho:\Hbf''\mono\Sp_{V''}$ along $F$ over $\Qbb$, where $\Hbf''$ is the spin $F$-group of $W''$ and $V''$ is an $F$-symplectic space such that $V=\Res_{F/\Qbb}V''$ and the base change of $\rho$ along $\sigma:F\mono\Rbb$ decomposes into a direct sum of copies of the unique spinor representation $P$ for $\Spin(2n+1,2)$.

We have explained that the restriction of $P$ to $\Spin(2n,2)$ is the direct sum of the two half-spin representations $P_-$ and $P_+$
of equal multiplicity, hence it suffices to apply \autoref{proposition restriction on the real part}\,(1)
to the further inclusion $\SU(n,1)\subset\Spin(2n,2)\subset\Spin(2n+1,2)$
to deduce the bound $d=[F:\Qbb]> 3+2^{-n-1}\big(\frac{1}{r}-\frac{1}{2}\big)$ from Xiao's inequality.

It remains to make precise the multiplicity $r$ along the parity of $N$. Assume that the inclusion $M'\mono\Acal_V$ is primary of multiplicity $m$: in other words, the rational symplectic representation defining $M'\mono\Acal_V$ is the scalar restriction of some $F$-linear symplectic representation $\Hbf\ra\Sp_W$, and $W\otimes_{F,\sigma}\Rbb$ is isomorphic to $m$-copies of the spinor representation of $\Spin(N,2)$. \begin{itemize}
\item If $N=1+2N'>2n$ is odd, with $N'\geq n$, then there exists only one spinor representation $P'$ for $\Hbf'(\Rbb,\sigma)\isom\Spin(1+2N',2)$, and its restriction to $\Hbf''(\Rbb,\sigma)\isom\Spin(2n,2)$ is isomorphic to $2^{N'-n}$ copies of $P_+\oplus P_-$. Since $P'$ is of multiplicity $m$ in $W\otimes_{F,\sigma}\Rbb$, this gives $r=m\cdot 2^{N'-n}$ for the multiplicity of $P_+\oplus P_-$ in the representation of $\Hbf''(\Rbb,\sigma)$.

\item If $N=2N'>2n$ is even, with $N'>n$, then either of the two half-spin representations $P'_+$ and $P'_-$ of $\Spin(2N',2)$ restricts to the unique spinor representation of $\Spin(2N'-1,2)$ and further to $2^{N'-n-1}$ copies of $P_+\oplus P_-$ of $\Spin(2n,2)$, and one obtains $r=m\cdot 2^{N'-n-1}$ for the multiplicity of $P_+\oplus P_-$ in the representation of $\Hbf''(\Rbb,\sigma)$.
\end{itemize}

Hence the multiplicity is always $r=m\cdot 2^{\lfloor(N-1)/2\rfloor-n}=m\cdot 2^{\lfloor(N+1)/2\rfloor-n-1}$, with $m$ the multiplicity of $W$, and the proof is completed.
\end{proof}

\vspace{5mm}
\section*{Acknowledgement}
It is our pleasure and honor to dedicate this work to Prof. Ngaiming Mok at his sixtieth anniversary. Prof. Mok has been well-known for his contribution to complex differential geometry and algebraic geometry, and his works have seen growing influences on young geometers. We congratulate him sincerely at this occasion, and wish him a fruitful career yet to come. %, just as the old Chinese saying: \begin{center}
%\emph{An steed in the stable still aspires to gallop a thousand li,}

%\emph{An aged hero retains his high aspiration even in high age.}
%\end{center}


\begin{thebibliography}{00}

\bibitem{chen lu zuo compositio}K. Chen, X. Lu, and K. Zuo, On the Oort conjecture for Shimura varieties of unitary and orthogonal types, Compositio Mathematica 152(2016), 889-917

\bibitem{deligne pspm} P. Deligne, Vari\'et\'es de Shimura, interpretation modulaire et construction de mod\`eles canoniques, Automorphic Forms, Representations, and L-functions, Part 2, Proceedings of Symposia in Pure Mathematics, vol. 33 (AMS Providence RI, 1979), 247-289

%\bibitem{fiori shimura subvarieties} A. Fiori, Shimura subvarieties of O(2,n)-type, preprint

\bibitem{hain locally symmetric} R. Hain, Locally symmetric families of curves and Jacobians, Aspects of Mathematics, vol. E33 (Vieweg, Braunschwig, 1999), 91-108

\bibitem{de jong zhang hilbert} J. de Jong, S. Zhang, Generic Abelian varieties with real multiplication are not Jacobians, Diophantine Geometry, CRM Series, vol. 4 (Ediziono della Normale, Pisa, 2007), 165-172

\bibitem{lu zuo} X. Lu and K. Zuo, The Oort conjecture for Shimura curves in the Torelli locus of curves, arXiv:1405:4751

\bibitem{moeller viehweg zuo} M. M\"oller, E. Viehweg, and K. Zuo, Stability of Hodge bundles and characterization of Shimura curves, Journal of Differential Geometry 92(2012), 71-151

\bibitem{mok rigidity book} N. Mok, Metric rigidity theorems on Hermitian locally symmetric manifolds, Series in pure mathematics, vol. 6, World Scientific

\bibitem{moonen oort survey} B. Moonen and F. Oort, The Torelli locus and special subvarieties, Handbook of Moduli, Vol. II, 549-594

%\bibitem{mueller-stach viehweg zuo} S. M\"uller-Stach, E. Viehweg, and K. Zuo, Relative proportionality for subvarieties of moduli spaces of K3-surfaces and abelian surfaces, Pure and Applied Mathematics Quarterly, 5(2009), no.3, 1161-1199

\bibitem{satake real} I. Satake, Holomorphic embeddings of symplectic domains into a Siegel space, Americal Journal of Mathematics, 87(1965), 425-461

\bibitem{satake rational} I. Satake, Symplectic representations of algebraic groups satisfying a certain analyticity condition, Acta Mathematica 117(1965), 215-279

%\bibitem{snow homogeneous} D. Snow, Homogeneous vector bundles, in Group Actions and Invariant Theory, CMS conference proceedings 10, 193-205

\bibitem{tsimerman andre oort} J. Tsimerman, A proof of the Andr\'e-Oort conjecture, cf. arXiv: 1506.01466

\bibitem{xiao slope} G. Xiao, Fibred algebraic surfaces with low slope, Mathematische Annalen, 276(1987), no.3, 449-466

\end{thebibliography}
\end{document}